\input epsf
\magnification=\magstep1


\hsize=6.35truein \vsize=9truein

\font\title=cmssdc10 at 15pt

\def\ee{\epsilon}

\def\frame #1{\vbox{\hrule height.1pt
\hbox{\vrule width.1pt\kern 10pt
\vbox{\kern 10pt
\vbox{\hsize 3.5in\noindent #1}
\kern 10pt}
\kern 10pt\vrule width.1pt}
\hrule height0pt depth.1pt}}
\def\bframe #1{\vbox{\hrule height.1pt
\hbox{\vrule width.1pt\kern 10pt
\vbox{\kern 10pt
\vbox{\hsize 5in\noindent #1}
\kern 10pt}
\kern 10pt\vrule width.1pt}
\hrule height0pt depth.1pt}}
\baselineskip 15pt

\centerline{\title A stochastic model of evolution}
\medskip
\centerline{Herv\'e Guiol, F\'abio P. Machado and Rinaldo B. Schinazi}
\centerline{TIMB-TIMC Univ. Grenoble, France, IME-USP, Brasil, Math. Dept. UCCS, USA}
\centerline{June 14 2010}
\bigskip
\footnote{}{Key words and phrases: Bak-Sneppen model, evolution, stochastic model, cut-off phenomenon}

\bigskip
\bigskip
{\bf Abstract.} We propose a stochastic model for evolution. Births and deaths of species occur
with constant probabilities. Each new species is associated with a fitness sampled from the uniform distribution on $[0,1]$.  Every time there is a death event then the type that is killed is the one with the smallest fitness.
We show that there is a sharp phase transition when the birth probability is larger than
the death probability. The set of species with fitness higher than a certain critical value
approach an uniform distribution. On the other hand all the species with fitness less than the
critical disappear after a finite (random) time.

\bigskip
{\bf 1. Introduction.}

Consider a discrete time model that starts from the empty set. At each time $n\geq 1$ with
probability $p$ there is a birth of a new species and with probability $q=1-p$ there is a death
of a species (if the system is not empty). Hence, the total number of species at time
$n$ is a random walk on the positive integers which jumps to the right with probability $p$ and
to the left with probability $q$. When the random walk is at 0 then it jumps to 1 with probability
$p$ or stays at $0$ with probability $1-p$.
 Each new species is associated with a random number. This random number is sampled from the uniform distribution on $[0,1]$. We think of the random number associated with a given species as being the fitness of the species. These random numbers are independent of each other and of everything else.
 Every time there is a death event then the type that is killed is the one with the smallest fitness.
 This is similar to a model introduced by Liggett and Schinazi (2009) for a different question.

 Take $p$ in $(1/2,1)$ and let
 $$f_c={1-p\over p}.$$
 Note that $f_c$ is in $(0,1)$.
 Let $L_n$ and $R_n$ be the set of species alive at time $n$ whose fitness is lower and higher than
 $f_c$, respectively. Since each fitness appears at most once almost surely we can identify each species to its fitness and think of $L_n$ and $R_n$ as sets of points in $(0,f_c)$ and $(f_c,1)$, respectively. Let $|A|$ denote the cardinal of set $A$. We are now ready to state our main result.
 \medskip
 {\bf Theorem. }{\sl Assume that $p>1/2$. Let $f_c={1-p\over p}.$

 (a)  The number $|L_n|$ of species whose fitness is below $f_c$ is a null recurrent birth and death chain. In particular, the set $L_n$ is empty infinitely often with probability one.

 (b) Let $f_c<a<b<1$ then
 $$\lim_{n\to\infty}{1\over n} |R_n\cap (a,b)|=p(b-a)\hbox { a.s.}$$
 }
 \medskip
 In words, there is a sharp transition at fitness $f_c$. No species with fitness below $f_c$ can survive forever. On the other hand species are asymptotically uniformly distributed on $(f_c,1)$.
 \medskip
\centerline{
\epsfxsize=6cm\epsfbox{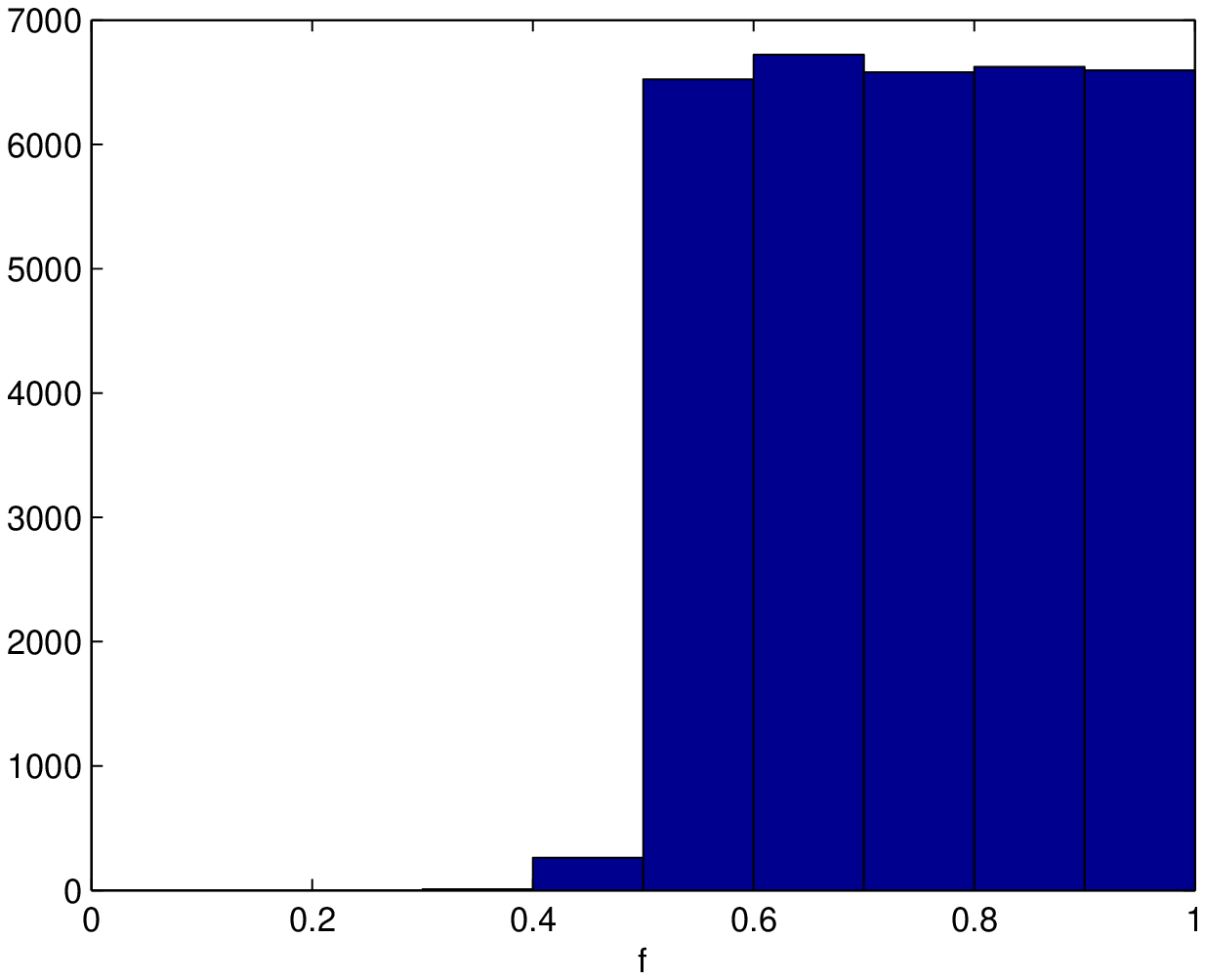}}
{\bf Figure 1.} This is the histogram of the fitnesses after 100,000 births and deaths for $p=2/3$.
We have $f_c=1/2$ and
as predicted by the Theorem the distribution on $(f_c,1)$ approaches an uniform.
\medskip
Observe that the larger $p$ is the more welcoming the environment is to new species. If $p$ is only
slightly larger than 1/2 then $f_c$ is close to 1 and only species with high fitness will survive.
On the other hand if $p$ is close to 1 then $f_c$ is close to 0 and even species with relatively low fitness will survive.

The 'kill the least fit' rule was introduced in the Bak-Sneppen model, see Bak and Sneppen (1993). In that model there is a fixed number $N$ of species
arranged in a circular graph. At each discrete time the site on the circle with the lowest
fitness and its two nearest neighbors have their fitness replaced by a random number independently sampled from the uniform distribution on $[0,1]$. Compared to our model there are two important differences: the number of species is fixed in the Bak-Sneppen model (in our model it is random) and there is some local interaction (kill the neighbors of the least fit). However, the same type of uniform behavior on some $(f_c,1)$ is expected for the Bak-Sneppen model but this is still unproved, see Meester and Znamenski (2003) and (2004).

In fact our result is more general than stated. The reader can easily check in our proof that there is nothing special about the uniform distribution. If fitnesses are sampled independently from the same fixed distribution then the limit in Theorem (b) is a.s. $pP(a<X<b)$ where $X$ is a random variable with the fixed fitness distribution. Based on computer simulations we conjecture that the same is true for the Bak-Sneppen model. There too the uniform distribution appears only because fitnesses are sampled from it.

\medskip
{\bf 2. Proof of the Theorem.}
\medskip
Part (a) is a well-known result for birth and death chains.
Recall that $L_n$ is the set of species whose fitness is lower than $f_c$ at time $n$. Observe that
 $|L_n|$ (the cardinal of $L_n$) increases by 1 with probability $pf_c$, decreases by 1 with probability
 $q$ (if it is not already at 0) and stays put with probability $p(1-f_c)$. Since $pf_c=q$,  it is easy to check that $|L_n|$ is null recurrent. See for instance Proposition II.2.4 in Schinazi (1999).
 \medskip
 We now turn to the proof of (b).
 Let $t_n$ be the number times $k\leq n$  for which $L_k$ is empty. That is,
$$t_n=|\{1\leq k\leq n:L_k=\emptyset\}|.$$
We will show that, for any $\ee>0$,  $t_n$ is almost surely less $n^{1/2+\ee}$ for $n$ large enough. The main step in the proof is the following Lemma.
\medskip
 {\bf Lemma. }{\sl There are positive constants $\gamma$ and $D$ such that for every $\ee>0$ we have
 $$P\Big(t_n>{2\over pf_c}n^{1/2+\ee}\Big)\leq D\exp(-\gamma n^{\ee}).$$
 }
 \medskip
{\it Proof of the Lemma.}

 Recall that we start from the empty set. After a geometric random time with mean ${1\over pf_c}$, denoted by $G_0$, the first species appears in $(0,f_c)$.
That is,
$$G_0=\min\{k\geq 1: L_k\not=\emptyset\}.$$
 Let
 $$E_1=\min\{k\geq G_0:L_k=\emptyset\}.$$
Hence, $E_1$ is the time it takes starting at time $G_0$ for $|L_.|$ to return to 0. More generally, we define for $i\geq 1$
$$G_i=\min\{k\geq G_0+E_1+\dots+G_{i-1}+E_i: L_k\not=\emptyset\},$$
and
$$E_{i+1}= \min\{k\geq G_0+E_1+\dots+E_i+G_i:L_k=\emptyset\}.$$

Note that the $(G_i)_{i\geq 0}$ and the $(E_i)_{i\geq 1}$ are two i.i.d. sequences. Moreover, the $G_i$ follow a geometric distribution with mean ${1\over pf_c}$.

Let $k_n$ be the number of times that $L_k$ hits the empty set by time $n$:
$$k_n=|\{2\leq k\leq n: |L_{k-1}|=1 \hbox{ and }|L_k|=0\}|.$$
That is, $k_n$ counts the number of times $L_k$ goes from 1 to 0 species for $k\leq n$.  Note that if $k_n=0$ then $t_n\leq G_0$.
Let $C={2\over pf_c}$. We now compute
 $$P(t_n>Cn^{1/2+\ee})\leq P(t_n>Cn^{1/2+\ee};k_n<n^{1/2+\ee})+P(k_n\geq n^{1/2+\ee}).\leqno (1)$$
For $k_n\geq 1$ we have
$$G_0+G_1+\dots+G_{k_n-1}<t_n\leq G_0+G_1+\dots+G_{k_n},$$
and for $k_n=0$ we have  $t_n\leq G_0$.
Hence,
$$P(t_n>Cn^{1/2+\ee};k_n<n^{1/2+\ee})\leq P(G_0+G_1+\dots+G_{m_n}>Cn^{1/2+\ee}),$$
where $m_n$ is the integer part of $n^{1/2+\ee}$. Now, the expected value of $G_0+G_1+\dots+G_{m_n}$ is ${m_n+1\over pf_c}$. By a large deviations inequality (see for instance Lemma (9.4) in Chapter 1 of Durrett (1996)) there exists $\gamma>0$ such that
$$P(G_0+G_1+\dots+G_{m_n}>Cn^{1/2+\ee})\leq \exp(-\gamma m_n)\leq \exp(-\gamma (n^{1/2+\ee}-1)).\leqno(2)$$
We now take care of the second term in the r.h.s. of (1).
Using that the $E_i$ are i.i.d. and that for $1\leq i\leq k_n-1$ they all must be less than $n$,
$$P(k_n\geq n^{1/2+\ee})\leq P(E_1<n)^{m_n-1}.$$
In order to estimate $P(E_1<n)$ we will compare $|L_n|$ to a simple symmetric random walk $W_n$ (one that jumps +1 or -1 with probability 1/2 at each step). We construct $W_n$ from $|L_n|$ by erasing the steps where $|L_n|$ stays put. If for instance we have $|L_1|=|L_2|=0$, $|L_3|=|L_4|=1$ and $|L_5|=2$ then we define $W_1=0$, $W_2=1$ and $W_2=2$. By doing so we get a simple symmetric random walk that visits the same sites (in the same order) as $|L_n|$  but in less time.
Hence, $|L_n|$ takes more time to go from 1 to 0 than $W_n$ does. Let $T_0$ be the time for $W_n$ to hit 0. We have
 $$P(k_n\geq n^{1/2+\ee})\leq P(E_1<n)^{m_n-1}\leq P_1(T_0<n)^{m_n-1}.$$
 It is well known that $P_1(T_0\geq n)$ is asymptotically $1/\sqrt{\pi n/2}$, see for instance Chapter III in Feller (1968). Hence, there are constants $\gamma'>0$ and $D$ such that
 $$P(k_n\geq n^{1/2+\ee})\leq \exp(-\gamma' {m_n-1\over n^{1/2}})\leq D\exp(-\gamma' n^\ee).\leqno(3)$$
Using (2) and (3) in (1) completes the proof of the Lemma.
\medskip
We are now ready to complete the proof of part (b) of the Theorem. Let $N_n$ be the total number of births in the model up to time $n$. Clearly, $N_n$ has a binomial distribution with parameters $n$ and $p$. Let $f_c<a<b<1$ we have
$$|R_n\cap (a,b)|\leq \sum_{i=1}^{N_n} 1_{(a,b)}(U_i),$$
where $1_{(a,b)}$ is the indicator function of the set $(a,b)$ and $(U_i)_{i\geq 1}$ is the sequence of fitnesses associated with births. Recall that the $U_i$ are i.i.d. and uniformly distributed on $(0,1)$. All this inequality is saying is that the number of points in the set $(a,b)$ at time $n$ is less than the number of births that occurred  up to time $n$ in the same set.

We now bound the number of deaths. We claim that the number of deaths in $(f_c,1)$ is at most $t_n$. This so because there can be a death in $(f_c,1)$ only when $(0,f_c)$ is empty and $t_n$ counts the number of times this happens up to time $n$. Hence,
$$\sum_{i=1}^{N_n} 1_{(a,b)}(U_i)-t_n\leq |R_n\cap (a,b)|\leq \sum_{i=1}^{N_n} 1_{(a,b)}(U_i).\leqno(4)$$
By the Law of Large Numbers,
$${1\over n}\sum_{i=1}^{N_n} 1_{(a,b)}(U_i)={N_n\over n}{1\over N_n}\sum_{i=1}^{N_n} 1_{(a,b)}(U_i)$$
converges a.s. to $pE(1_{(a,b)}(U))=p(b-a).$

 On the other hand, by our Lemma and  the Borel-Cantelli Lemma  there is  almost surely a natural $N$ such that $t_n\leq {2\over pf_c}n^{1/2+\ee}$ for $n\geq N$. In particular, $t_n/n$ converges to 0 a.s. We use the two preceding limits in (4) to conclude that a.s.
 $$\lim_{n\to\infty}{1\over n}|R_n\cap (a,b)|=p(b-a).$$
 This completes the proof of the Theorem.
 \medskip
 {\bf Note} We have just learned that Ben-Ari, Matzavinos and Roitershtein (2010) proved a central limit
theorem and a law of the iterated logarithm for our model, by
developing further some ideas presented in this paper.

\medskip
{\bf Acknowledgements} H.G. thanks IXXI for partial support, Math. Department of University of
Colorado at Colorado Springs and IME USP Brasil for kind hospitality. F.M. thanks Math. Department
of University of Colorado at Colorado Springs for partial support and kind hospitality. R.B.S. was partially supported by N.S.F. grant DMS-0701396.

\bigbreak
{\bf References.}

P. Bak and K. Sneppen (1993). Punctuated equilibrium and criticality in a simple model of evolution. {\sl Phys. Rev. Lett.}, {\bf 74}, 4083-4086.

I. Ben-Ari, A. Matzavinos and  A. Roitershtein (2010). On a species survival model.
arXiv:1006.2585

R. Durrett (1996). {\sl Probability: theory and examples} (second edition). Duxbury Press.

W. Feller (1968). {\sl An introduction to probability theory and its applications}. Volume 1. John Wiley (3rd edition).

T. Liggett and R. B. Schinazi (2009). A stochastic model for phylogenetic trees. {\sl J. Appl. Prob.}, {\bf 46}, 601-607.

R. Meester and D. Znamenski (2003). Limit behavior of the Bak-Sneppen evolution model. {\sl Ann. Prob.}, {\bf 31}, 1986-2002.

R. Meester and D. Znamenski (2004).
Critical thresholds and the limit distribution in the Bak-Sneppen model.
{\sl Com. Math.Phys.} {\bf 246}, 63-86.

R.B.Schinazi (1999). {\sl Classical and spatial stochastic processes.}  Birkhauser.

\end